\documentclass{jnmp}

\usepackage{amsmath}
\usepackage{graphicx}
\usepackage{epsfig}

\JNMPnumberwithin{equation}{section}

\theoremstyle{definition}

\begin{document}

\renewcommand{\evenhead}{I~G Korepanov and E~V~Martyushev}
\renewcommand{\oddhead}{Distinguishing Three-Dimensional Lens Spaces $L(7,1)$ and $L(7,2)$}

\thispagestyle{empty}

\FirstPageHead{9}{1}{2002}{\pageref{korepanov-firstpage}--\pageref{korepanov-lastpage}}{Article}

\copyrightnote{2002}{I~G Korepanov and E~V~Martyushev}

\Name{Distinguishing Three-Dimensional Lens Spaces {\mathversion{bold}$L(7,1)$}
 and {\mathversion{bold}$L(7,2)$}
by Means of Classical Pentagon Equation}

\label{korepanov-firstpage}

\Author{I~G KOREPANOV~$^*$ and E~V~MARTYUSHEV}

\Address{South Ural State University, 76 Lenin avenue,
454080 Chelyabinsk, Russia\\
$^*$~E-mail: kig@susu.ac.ru}

\Date{Received May 15, 2001; Revised October 20, 2001;
Accepted October 26, 2001}

\begin{abstract}
\noindent
We construct new topological invariants of three-dimensional manifolds
which can, in particular, distinguish homotopy equivalent lens spaces
$L(7,1)$ and $L(7,2)$. The invariants are built on the base of a classical
(not quantum) solution of pentagon equation, i.e.\ algebraic
relation corresponding to a ``2 tetrahedra $\to$ 3 tetrahedra'' local
re-building of a manifold triangulation. This solution, found earlier by
one of the authors, is expressed in terms of metric characteristics
of Euclidean tetrahedra.
\end{abstract}

\section{Introduction}

The present work is a direct continuation of paper~\cite{1}. We recall that
a new construction was proposed in it of an invariant of 3-dimensional
closed oriented piecewise-linear manifolds.
The invariant was built on the base of a classical
(not quantum) solution of pentagon equation, i.e.\ algebraic
relation corresponding to a ``2 tetrahedra $\to$ 3 tetrahedra'' local
re-building of a manifold triangulation
(formula (5) in~\cite{1}. Below, we are using notations
like~\cite[(5)]{1} for references of such kind).

Only the simplest realization of these ideas was done in paper~\cite{1}.
An invariant~$I$ was constructed that was expressed in terms of Euclidean
metric values assigned to every tetrahedron. Then $I$ was calculated for
the sphere~$S^3$, the projective space~$\mathbb RP^3$, and also for
some lens spaces $L(p,q)$ (it was not described in paper~\cite{1} how
the calculations were done for $L(p,q)$. It will suffice for us now if
we say that they were a direct generalization of calculations
for~$\mathbb RP^3$).

For all tested examples of lens spaces we got
\begin{gather}
I\bigl(L(p,q)\bigr)=\frac{1}{p^3}.
\tag{\cite[(39)]{1}}
\end{gather}
This showed that our invariant~$I$ in such form as presented in~\cite{1}
was not very interesting for topology because it could not
distinguish the lens spaces with different~$q$.
At the same time, in the concluding Section~7 of paper~\cite{1} an idea was
proposed how to modify~$I$ in the way that could hopefully
lead to a more refined
invariant. The idea consists, briefly, in using some mappings of the
manifold's {\em universal covering\/} in the Euclidean space~$\mathbb R^3$
in order to put metric values in correspondence to tetrahedra.

In the present work we make this idea concrete and show that in such way
one can indeed construct more refined invariants. Our key
example will be lens spaces $L(7,1)$ and $L(7,2)$.
As is well-known, they are homotopy equivalent but not homeomorphic,
and can serve as a natural test of nontriviality
for topological invariants. We would like to say at once that our new
invariants pass this test and thus can, most likely, compete successfully
with the quantum invariants of 3-dimensional manifolds.
A remark on possible advantages of classical invariants is made
in the concluding Section~\ref{sec discussion}.

We consider triangulations or ``pre-triangulations'' of 3-dimensional
closed oriented manifolds. By pre-triangulation we
understand a manifold represented as a ``pre-simplicial
complex'' in the sense of textbook~\cite{HW}: the manifold is decomposed
in simplices of dimensions 3, 2, 1 and~0 in such way that
the boundary of a given
simplex may consist of {\em coinciding\/} simplices of lower dimensions
(in particular, in Section~\ref{sec linzy} we decompose lens
spaces in tetrahedra in such way that only two vertices of each
tetrahedron are different and, moreover, there are sometimes
coinciding ones between the six edges of a tetrahedron. Moreover, the
simplices of dimensions~$>0$ may not
be determined by the set of their vertices).

We are going to put metric values (edge lengths, dihedral angles, volumes)
in correspondence to all tetrahedra in the pre-triangulation, as it was done
in paper~\cite{1}. To be exact, we will do the following.
Choose a consistent orientation of all tetrahedra in the
pre-triangulation (using the manifold's orientability). Assign to all
edges (one-dimensional cells of the complex) positive real numbers~--- lengths~---
in such way that every individual tetrahedron could be
placed in a 3-dimensional Euclidean space. This means that one can calculate
in a usual way dihedral angles at each edge in every tetrahedron,
and one can assume that those angles take values from
$0$ to~$\pi$. Next, we assign to each
tetrahedron a ``$+$'' or ``$-$'' sign and impose the following
condition: the algebraic sum of dihedral angles around
each edge must equal $0$ modulo~$2\pi$, where ``algebraic sum''
means that each dihedral angle, which is calculated, as we have
mentioned, from
the six edge lengths of some tetrahedron, is taken with the sign
``$+$'' or ``$-$'' that is assigned to that tetrahedron.
This is exactly the condition that is implied in formulas
of paper~\cite{1}, starting from~\cite[(5)]{1}.

Such assignment (of lengths to edges and signs to tetrahedra)
was realized in paper~\cite{1} in the following way. All
{\em vertices\/}  of the complex were mapped into the 3-dimensional Euclidean
space~$\mathbb R^3$. For every tetrahedron in the complex we took the
images of its four vertices in~$\mathbb R^3$ and demanded that the
tetrahedron spanned by them should be
nondegenerate (of course, this is possible
only if all the four vertices are different). Assume also that we have
fixed the orientation of~$\mathbb R^3$. Then a tetrahedron in the
pre-triangulation either conserves its orientation under the mapping
to~$\mathbb R^3$ or changes it (recall that we have fixed
a consistent orientation for the tetrahedra). In the first case, we assign
the ``$+$'' sign to it, in the second~--- the ``$-$'' sign.
As for the edge lengths, they are, of course,
the distances between the relevant vertices in~$\mathbb R^3$.

In the present paper, we generalize this construction according
to the remark in Section~7 of paper~\cite{1} (and it is exactly
this generalization that leads us to our new results). We take
the {\em universal covering\/} of the
pre-triangulation of manifold~$M$ considered as a cell complex.
Its vertices are divided in classes~--- the inverse images
of each given vertex with respect to the covering map.
Choose a representative in each class and map it into
some point in the space~$\mathbb R^3$. Next, fix a homomorphism
$\varphi\colon\; \pi_1(M)\to E_3$ of the fundamental group of manifold~$M$
in the group of motions of 3-dimensional Euclidean space. Suppose
that some vertices
$A_1$ and $A_2$ in the universal covering belong to the same class,
then $A_2=g(A_1)$, where $g\in \pi_1(M)$. Denote the images of vertices
in $\mathbb R^3$ by the same letters with the tilde. We demand
that for all such vertices $\tilde A_2=\varphi(g)\,(\tilde A_1)$.

This construction yields (in the general position) a one-to-one
correspondence between the vertices in the universal covering of our
complex and their images in~$\mathbb R^3$. As for simplices
of dimension~$>0$, each of them yields thus
a simplex in~$\mathbb R^3$~--- the convex hull of corresponding
vertices (we emphasize that we do not care about possible intersections
of such simplices in~$\mathbb R^3$ but we demand that all the vertices
of each tetrahedron in the {\em universal covering\/} should be different). Hence
we obtain lengths for the edges of the universal covering and signs~---
for its 3-dimensional cells in the same way as it was done
two paragraphs before. Then, it is clear that if two edges
are mapped in the same one by the covering map then their lengths
obtained from our construction are equal (because their images
in~$\mathbb R^3$ are taken one into another by elements $\varphi(g)\in E_3$).
The same holds for the tetrahedron signs. Thus, we have assigned correctly
the lengths to the edges and the signs to the tetrahedra of
the manifold~$M$'s
pre-triangulation itself and not only to its universal covering.

In Section~\ref{sec linzy} we present a concrete example of this
construction for the case of lens spaces.

It can be shown that our construction with the universal covering yields
all possible ways of assigning lengths to the edges and signs
to the tetrahedra with the
condition that the algebraic sum of dihedral angles around each edge
must be $0$ modulo~$2\pi$. We will not,
however, use this fact in the present paper.

{\bf Important remark:} A homomorphism
$
\varphi\colon\;\pi_1(M) \to E_3
$
takes part in our construction. This will lead to the fact that our
invariant corresponds in reality to a
{\em pair\/} $(M,\varphi)$. In particular, for $M=L(7,q)$ there exist,
as we will see, three non-equivalent $\varphi$, excluding the trivial
homomorphism whose image is unity (in~\cite{1} we were dealing,
in essense, with exactly this ``trivial'' case).

Thus, we have assigned metric values to elements of pre-triangulation
in a more general way than in paper~\cite{1}. Most of the other constructions
in~\cite{1} do not need any generalization: the main thing that we are
going to do is calculate our key differential form~\cite[(30)]{1} and
extract a numeric invariant from it.

The contents of the remaining sections is as follows. In
Section~\ref{sec obsch-proizv} we write down general formulas for the
partial derivatives of ``defect angles'' (we recall {\it en passant\/}
their definition) with respect to the edge lengths. The matrix of
such derivatives is needed for constructing our invariants. In
Section~\ref{sec linzy} we are concerned with calculations as such
for lens spaces $L(7,1)$ and~$L(7,2)$. In the concluding
Section~\ref{sec discussion} we put forward a conjecture
concerning the value of
our invariants for {\em any\/} lens spaces, remark on
the possible advantages of classical invariants,
and discuss some intriguing unsolved problems.

\section{Derivatives of defect angles w.r.t.\ edge lengths}
\label{sec obsch-proizv}

In the Introduction we mapped the universal covering of a pre-triangulated
manifold~$M$ in Euclidean space~$\mathbb R^3$. The pull-back of one cell
(simplex) from~$M$ with respect to the covering map was a set of
$N={\mathop{\rm card}\nolimits}\,\bigl(\pi_1(M)\bigr)$ cells which were mapped, according to our
construction, in isometric simplices in~$\mathbb R^3$.

Suppose that Euclidean tetrahedra
$\tilde A_1\tilde B_1\tilde C_1\tilde D_1$, $\ldots\,$,
$\tilde A_N\tilde B_N\tilde C_N\tilde D_N$ correspond in such way to
a tetrahedron $ABCD$ lying in~$M$.
Now we, first, will no longer write the
tildes (that stayed for points in $\mathbb R^3$, see Introduction. We thus
identify the vertices in the universal
covering and their images in~$\mathbb R^3$). Second, we define
the oriented volume of tetrahedron $ABCD$ as $1/6$ of the
following triple scalar product:
\begin{equation}
6V_{ABCD}=\overrightarrow{A_1B_1}\,\overrightarrow{A_1C_1}\,
\overrightarrow{A_1D_1}\ \left(=\cdots =\overrightarrow{A_NB_N}\,\overrightarrow{A_NC_N}
\,\overrightarrow{A_ND_N}\,\right).
\label{eq or-ob'em}
\end{equation}

Our next task is to make some formulas of paper~\cite{1} more accurate.
Namely, we are going to take carefully into account the orientation
of angles and volumes and thus remove the absolute value signs.
Introduce (small) defect angles aroung edges like in paper~\cite{1}:
slightly (and otherwise arbitrarily) change the edge lengths; after that
the algebraic sums of dihedral angles aroung edges cease to be zero.
By definition, such sum around a given edge~$a$ is $-\omega_a$,
where $\omega_a$ is called defect angle (defined modulo~$2\pi$).
Note that the sign of expression (\ref{eq or-ob'em}) is nothing
but the sign corresponding to the tetrahedron, i.e.\ the sign
with which the dihedral angles of the tetrahedron enter in
each algebraic sum~$(-\omega_a)$.

We will need formulas for the partial derivatives
$\partial \omega_a/\partial l_b$
of defect angles with respect to edge lengths taken when all $\omega_a=0$.
We are going to express these derivatives in terms of
those lengths and oriented tetrahedron volumes. Nonzero derivatives
are obtained in the following cases.

\subsubsection*{1st case}

Edges $a=DE$ and $b=AB$ are skew edges of tetrahedron $ABDE$,
and there is no more tetrahedron in the pre-triangulation that would
have as its edges both $a$ and~$b$:
\begin{equation}
\frac{\partial \omega_{DE}}{\partial l_{AB}} = -\frac {l_{AB}\,l_{DE}}{6}\,
\frac {1}{V_{ABED}}.
\label{eq poka 2}
\end{equation}

This formula coincides, in essense, with formula~\cite[(3)]{1}. The only
difference is that we now take into account the signs of angles
and volumes.

\begin{figure}[th]
\begin{minipage}[t]{0.46\linewidth}
\centerline{\epsfig{figure=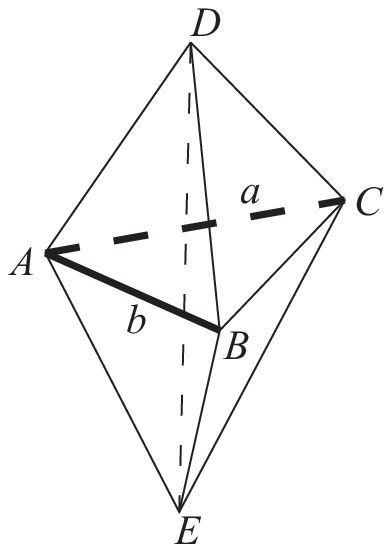}}
\caption{Illustration to formula~(\protect\ref{eq *_0}).}
\label{fig peresek rebra}
\end{minipage}\hfill
\begin{minipage}[t]{0.46\linewidth}
\centerline{\epsfig{figure=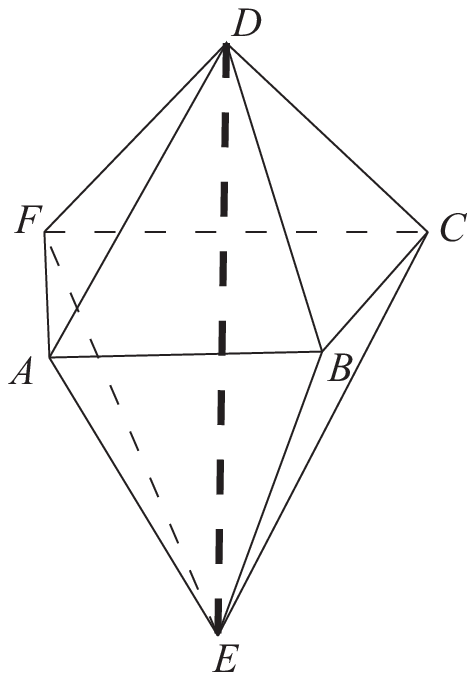}}
\caption{Illustration to formula~(\protect\ref{eq 4j}).}
\label{r2}
\end{minipage}
\end{figure}

\subsubsection*{2nd case}

Edges $a=AC$ and $b=AB$ belong to the common face of two neighboring
tetrahedra and, again, there is no more tetrahedron that would have both
$a$ and~$b$ as its edges (Fig.~\ref{fig peresek rebra}):
\begin{equation}
\frac{\partial \omega_{AC}}{\partial l_{AB}} = \frac {l_{AB}\,l_{AC}}{6}\,
\frac {V_{BCED}}{V_{ABCD}\,V_{EABC}}.
\label{eq *_0}
\end{equation}

To prove this formula, we draw edge~$DE$ in
Fig.~\ref{fig peresek rebra}. Suppose that lengths $l_a$ and $l_{DE}$
in Fig.~\ref{fig peresek rebra} are free to change,
while the lengths of remaining eight edges are fixed. We get:
\begin{equation}
\frac{\partial l_{DE}}{\partial l_b} = -\frac{l_b}{l_{DE}} \frac{V_{CEDA}\, V_{BCED}}{V_{ABCD}\,
V_{EABC}}
\label{eq *_1}
\end{equation}
(this is, in essense, formula~\cite[(1)]{1} but
now with the proper notice to volume signs).

Denote $\varphi_{AC}$ the dihedral angle at edge $AC$ in tetrahedron $CEDA$.
Then
\begin{equation}
\frac{\partial \varphi_{AC}}{\partial l_{DE}} = \frac{l_{AC}\, l_{DE}}{6V_{CEDA}}.
\label{eq *_2}
\end{equation}

Multiplying (\ref{eq *_1}) by (\ref{eq *_2}) and using the fact that
$\partial \omega_{AC}
/\partial l_b = - \partial \varphi_{AC} / \partial l_b$, we get
formula~(\ref{eq *_0}).

\subsubsection*{3rd case}

Edge $a=b=DE$ is common for exactly three tetrahedra
$ABED$, $BCED$ and $CAED$, and no one of those
tetrahedra contains this edge more than one time:
\begin{equation}
\frac{\partial \omega_{DE}}{\partial l_{DE}} = -\frac {l_{DE}^2}{6}\,
\frac {V_{ABCD}\,V_{EABC}}{V_{ABED}\,V_{BCED}\,V_{CAED}}.
\label{eq poka 6}
\end{equation}

This is simply formula~\cite[(5)]{1} but with all signs taken into account.

\subsubsection*{4th case}
Edge $a=b=DE$ is common for $>3$ tetrahedra,
and again no one of those
tetrahedra contains this edge more than one time.

We write out this formula only for the case of four tetrahedra (in the case
of their greater number, the generalization is obvious)
depicted in Fig.~\ref{r2}:
\begin{equation}
\frac{\partial \omega_{DE}}{\partial l_{DE}} = -\frac {l_{DE}^2}{6} \left(
\frac {V_{ABCD}\,V_{EABC}}{V_{ABED}\,V_{BCED}\,V_{CAED}}+
\frac {V_{ACFD}\,V_{EACF}}{V_{ACED}\,V_{CFED}\,V_{FAED}} \right).
\label{eq 4j}
\end{equation}
Here in the right-hand side we wrote the right-hand side of~(\ref{eq poka 6})
plus the similar term obtained from it by replacing
$C\to F$, $B\to C$. Formula~(\ref{eq 4j}) comes out if we
draw the diagonal $AC$ in Fig.~\ref{r2}  and apply
formula~(\ref{eq poka 6}) to each of figures
$ABCDE$ and~$ACFDE$. Adding up the defect angles around $DE$ in each of
those figures we get nothing else than $\omega_{DE}$ for the
whole figure $ABCFDE$,
because the ``redundant'' dihedral angle, namely the angle at~$DE$ in
tetrahedron $AECD$, enters two times with different signs.

\subsubsection*{5th case}

In reality, our definition of pre-triangulation admits such glueings
between faces and, consequently, edges of
tetrahedra where different
{\em combinations\/} of above mentioned cases appear. To be exact,
the edge~$b$ may
\vspace{-2mm}
\begin{itemize}
\itemsep0mm
\item
lie opposite to edge $a$ (1st case) in one or more than one
tetrahedron;
\item
belong to the same two-dimensional face as $a$ (2nd case), and again
there may be more than one such faces;
\item
coincide with $a$ (as in the 3rd or 4th case),
\vspace{-2mm}
\end{itemize}
and these possibilities (as we will see in Section~\ref{sec linzy})
do not exclude one another.
It is also worth mentioning that if a tetrahedron contains a given edge~$a$
several times then $(-\omega_a)$ includes, of course, the {\em sum
of all\/} corresponding
dihedral angles of that tetrahedron multiplied by the sign of its
oriented volume~(\ref{eq or-ob'em}).

Thus, there may be several ``ways of influence'' of differential $dl_b$
on differential~$d\omega_a$. Clearly, the terms corresponding to those ways
sum together. Formulas (\ref{eq pok1}) and~(\ref{eq pok2}) below are good
examples of how the expression for $\partial \omega_a / \partial l_b$
(in them, it is divided by~$l_al_b$) can look like.

The partial derivatives that we were dealing with in this section
will be used to calculate the
differential form \cite[formula (30)]{1} from which we will extract
manifold invariants. The formulas and arguments of paper~\cite{1} leading to that
formula remain valid if we only modify the notion of
``permitted length configuration'', extending the word ``permitted''
to any set of edge lengths obtained according to the construction
with universal covering presented in
the Introduction to this work.

The remarks in the end of Section~4 of~\cite{1} (after Theorem~4) on the
behavior of form \cite[(30)]{1} under the adding of a new vertex to the
(pre-)triangulation remain valid as well. In formulas
\cite[(31) and (32)]{1}, one can take as point~$E$ any one
of ``copies'' $E_1, \ldots, E_N$~--- images in $\mathbb R^3$ of its
pull-back with respect to the covering map. On the other hand, we will
see in the next section that we must change the ``standard'' differential
form \cite[(33)]{1} in order to calculate our new
invariants for $L(7,1)$ and~$L(7,2)$.

\section{Calculation of invariants for {\mathversion{bold}$L(7,1)$} and {\mathversion{bold}$L(7,2)$}}
\label{sec linzy}

\begin{figure}[th]
\centerline{\epsfig{figure=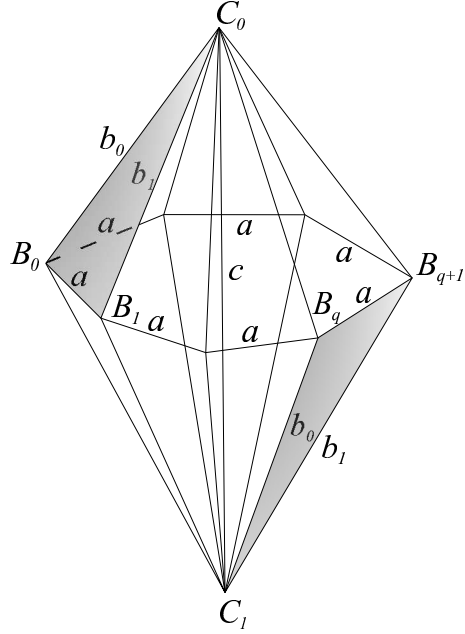}} \caption{Lens
space as a bipyramid with identified faces.} \label{ris3}
\end{figure}

We recall some facts about lens spaces.
Let $p > q > 0$, $p \geqslant 3$, be a pair of relatively prime integers.
Consider a $p$-gonal bipyramid, i.e.\ the union of two cones over
a regular $p$-gon. Denote as $B_0, B_1,\ldots, B_{p-1}$ the vertices
of the $p$-gon, and by $C_0$ and $C_1$~--- the cone vertices.
For every $i$, we glue face $B_iC_0B_{i+1}$ to face $B_{i+q}C_1B_{i+q+1}$
(the subscripts are taken modulo $p$; the vertices are glued in the
order in which they are written). What we thus obtain is exactly the
lens space $L(p,q)$ (Fig.~\ref{ris3}).

$L(p,q)$ is a 3-dimensional closed oriented manifold.
Its fundamental group $\pi_1\bigl(L(p,q)\bigr)$ is isomorphic to $\mathbb{Z}_p$.
As a generator of $\pi_1\bigl(L(p,q)\bigr)$ one can take a path going
along the edge $a=B_iB_{i+1}$ (note that points $B_i$ as well as edges
$B_iB_{i+1}$ are identified among themselves for all~$i$), or along
the edge~$c=C_0C_1$.

The spaces $L(7,1)$ and $L(7,2)$ that we are going to consider are
interesting because they are homotopy equivalent~\cite{HW} but not
homeomorphic~\cite{milnor}. We will obtain one possible proof of their
non-homeomorphness when we find out that our invariants
for those spaces do not coincide.

A homomorphism of the fundamental group of a lens space
in the group~$E_3$ of motions of Euclidean space can be characterized
by an integer~$k$ if we demand that $a\in \pi_1\bigl(L(p,q)\bigr)$
be mapped in the rotation through angle $2\pi k/p$ around some axis.
That axis can be taken by definition as the coordinate axis~$z$.
There exist three essentially different nonzero $k$ for $p=7$ ($k=0$
leads to the case of paper~\cite{1}): $k=1,2$ and~$3$.

The image of the universal covering of space $L(p,q)$ (if $k\ne 0$) in
$\mathbb R^3$ looks as follows. Vertices $B_j$ (``copies'' of one and
the same vertex~$B$) are situated in vertices of a regular $p$-gon
for which axis~$z$ is the symmetry axis. The same applies to vertices
$C_j$ (``copies'' of vertex~$C$. Thus, our bipyramid has only
two different [not identified with each other] vertices: $B$ and~$C$).
We can assume that the cylindrical coordinates of points $B_j$ and $C_j$
are:
\begin{equation}
B_j \left( \rho, \frac{2\pi jk}{p},0 \right),\qquad
C_j \left( \sigma, \alpha+\frac{2\pi qjk}{p},s \right),
\label{eq *_10}
\end{equation}
where $\rho>0$, $\sigma>0$, while $\alpha$ and $s$ are any real numbers.

Return to Fig.~\ref{ris3}. Join points $C_0$ and $C_1$ with an edge.
Then the bipyramid finds itself decomposed in $p$ tetrahedra
grouped around edge~$C_0C_1$.
Assume the following notations for edges (notation ``$a$'' and ``$c$''
have been already introduced but we repeat them):
\begin{gather}
a = B_0 B_1 = \cdots = B_{p-1} B_0, \\
b_j = C_0 B_j = C_1 B_{j+q \!\!\!\! \pmod{p}}, \quad j=0,\dots, p-1, \\
c = C_0 C_1.
\end{gather}
On the whole, there are $p+2$ edges, and this will be the size of square
matrix
$A=(\partial\omega_i/\partial l_j)$ entering in formula~\cite[(30)]{1}.
Both rows and columns of this matrix correspond to edges
$a,b_0,\ldots,b_{p-1},c$.
Decompose this set of edges in subsets ${\cal C}=\{b_1,\ldots,b_{p-2}\}$
and $\overline{\cal C}=\{a,b_0,b_{p-1},c\}$. Recall our description of the image
of universal covering of a
lens space in~$\mathbb R^3$. It is obvious from geometric considerations
that this image is determined by four parameters, say,
$\rho$, $\sigma$, $s$ and $\alpha$ in formulas (\ref{eq *_10}), or otherwise
by the lengths of four edges which we have
included in~$\overline{\cal C}$.
The fact important for us is that these four lengths (in the general position) can take
arbitrary infinitesimal increments from which
the other length differentials are determined unambiguously,
and one can assume that they are determined from conditions $d\omega_i=0$,
where the letter~$i$ runs over all edges. According to paper~\cite{1},
this means that our decomposition of the set of edges in $\cal C$
and $\overline{\cal C}$ is suitable for use in formula~\cite[(30)]{1}.

We will slightly modify that formula by introducing a matrix
\begin{equation}
F\stackrel{\rm def}{=} \left(\frac{1}{l_il_j} \frac{\partial\omega_i}{\partial l_j}\right)
\end{equation}
instead of $A$. Then the form~\cite[(30)]{1} is written as
\begin{equation}
\left| \frac{\displaystyle\bigwedge_{\overline{\cal C}} l\, dl}{\,\root\of
{\left| \det \left( F|_{\cal C} \right) \prod 6V
\vrule height 1.8ex depth 0.9ex width 0pt
\right|}} \right|,
\label{eq new30}
\end{equation}
where, of course, $F|_{\cal C}$ is the diagonal submatrix of
matrix~$F$ corresponding to subset~$\cal C$.

Our next task is to calculate explicitly the matrix elements of
$F|_{\cal C}$ for $L(7,1)$ and $L(7,2)$.
We will explain how the calculations are done on the example of two
types of matrix elements of $L(7,1)$, while
we will simply write out the result for all other matrix elements.

At the moment, we are considering the case of arbitrary $p$ and $q$.

Introduce a brief notation
\begin{equation}
V_i = V_{C_0 C_1 B_{i+1} B_i}.
\label{eq opredelenie V_i}
\end{equation}
An easy calculation shows that
\begin{equation}
V_i = \frac{1}{6}\, R\, \sin\left( \alpha + \frac{\pi k(q-1-2i)}
{p}\right),
\label{eq V_i}
\end{equation}
where we have denoted
\begin{equation}
R = 4 \rho \sigma s\,\sin \frac{\pi k}{p}\,\sin \frac{\pi qk}{p}.
\label{eq *_R}
\end{equation}

Now set $p=7$, $q=1$.

\begin{figure}[th]
\begin{minipage}[t]{0.46\linewidth}
\centerline{\epsfig{figure=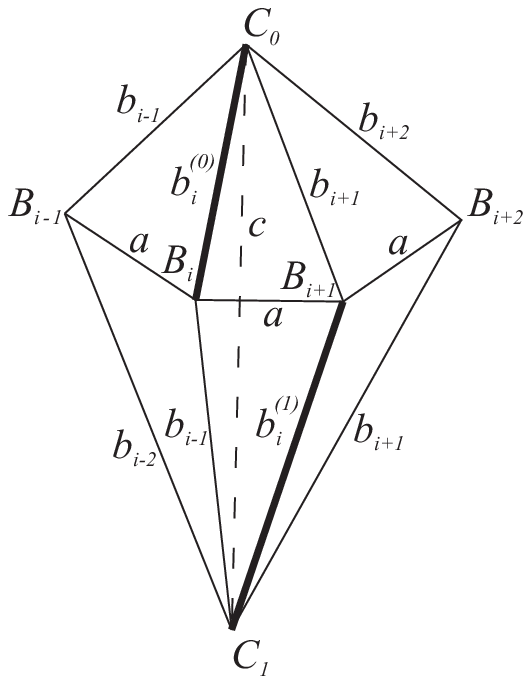}} \caption{A
fragment of bipyramid~--- illus\-tration to
formula~(\protect\ref{eq pok1}). Two copies of edge~$b_i$ are
drawn in thick lines.} \label{1-j kusok bipiramidy}
\end{minipage}\hfill
\begin{minipage}[t]{0.46\linewidth}
\centerline{\epsfig{figure=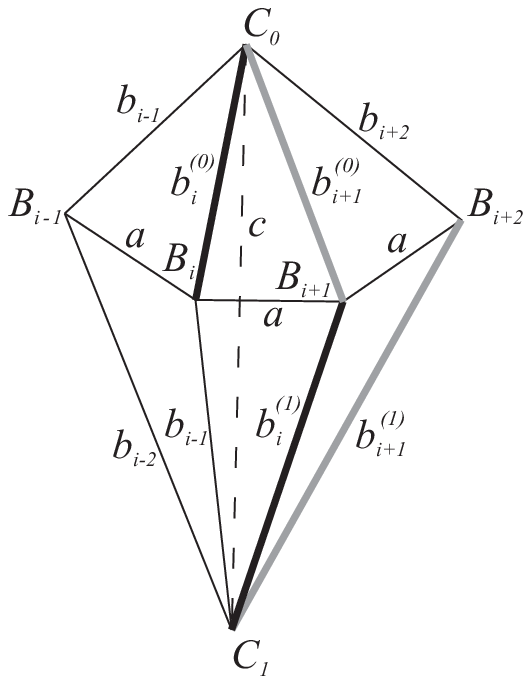}} \caption{A
fragment of bipyramid~--- illus\-tra\-tion to
formula~(\protect\ref{eq pok2}). Different types of thick lines
show edges $b_i$ and~$b_{i+1}$ (each in two copies glued
together).} \label{2-j kusok bipiramidy}
\end{minipage}
\end{figure}

\subsubsection*{Diagonal matrix elements}
\[
\left(F|_{\cal C}\right)_{i,i}=
\frac{1}{l_{b_i}^2}\frac{\partial \omega_{b_i}}{\partial l_{b_i}},\qquad i=1,\ldots,5.
\]

In Fig.~\ref{1-j kusok bipiramidy}
a fragment of bipyramid is presented containing both copies of edge~$b_i$
that are glued together.
We denote them $b_i^{(0)}$ and~$b_i^{(1)}$. Each of them belongs to two
tetrahedra, thus, the total number of tetrahedra around edge~$b_i$ is
four. This would be the
4th case from Section~\ref{sec obsch-proizv} be it not combined
with the 1st case (thus giving the 5th case as a result), because edge
$b_i$ turns out also to lie opposite itself
in tetrahedron $C_0C_1B_{i+1}B_i$. This leads to the adding of two more identical
terms of type~(\ref{eq poka 2}) (the fact that there are two such
terms can be explained as follows:
$b_i^{(0)}$ influences the dihedral angle at~$b_i^{(1)}$, while $b_i^{(1)}$
influences the dihedral angle at~$b_i^{(0)}$).

The result is the following combination of formulas (\ref{eq poka 6})
and~(\ref{eq poka 2}):
\begin{gather}
\left(F|_{\cal C}\right)_{i,i}=-\frac{1}{6}\left(
\frac{V_{C_0C_1B_{i+2}B_i} V_{C_0B_iB_{i+1}B_{i+2}}}
{V_{C_0C_1B_{i+2}B_{i+1}} V_{C_0C_1B_{i+1}B_i} V_{C_1B_iB_{i+1}B_{i+2}}}
\right. \nonumber\\
\phantom{\left(F|_{\cal C}\right)_{i,i}=}
\left. {}+ \frac{V_{C_0C_1B_{i+1}B_{i-1}} V_{C_1B_{i-1}B_iB_{i+1}}}
{V_{C_0C_1B_iB_{i-1}} V_{C_0C_1B_{i+1}B_i} V_{C_0B_{i-1}B_iB_{i+1}}} +
\frac{2}{V_{C_0C_1B_{i+1}B_i}} \right).
\label{eq pok1}
\end{gather}

Formula (\ref{eq pok1}) can be simplified considerably
if we take into account that points $B_j$ and $C_j$ are
located in two parallel planes. This gives, firstly:
\[
V_{C_0B_iB_{i+1}B_{i+2}}=V_{C_1B_iB_{i+1}B_{i+2}},\qquad
V_{C_0B_{i-1}B_iB_{i+1}}=V_{C_1B_{i-1}B_iB_{i+1}},
\]
which leads at once to partial cancellation of numerators and
denominators in two first terms in parentheses in the
right-hand side of~(\ref{eq pok1}). Further simplification of those
very terms is achieved by noting that the volumes that remain
in their numerators are expressed through the volumes in their
denominators as follows:
\begin{gather*}
V_{C_0C_1B_{i+2}B_i}=V_{C_0C_1B_{i+1}B_i}+V_{C_0C_1B_{i+2}B_{i+1}},\\
V_{C_0C_1B_{i+1}B_{i-1}}=V_{C_0C_1B_iB_{i-1}}+V_{C_0C_1B_{i+1}B_i}.
\end{gather*}
Finally we get
\begin{equation}
\left(F|_{\cal C}\right)_{i,i}=-\frac{4}{6V_i} -\frac{1}{6V_{i-1}}
-\frac{1}{6V_{i+1}}.
\end{equation}

\subsubsection*{Next-to-diagonal matrix elements}
\[
\left(F|_{\cal C}\right)_{i,i+1}=
\frac{1}{l_{b_i}l_{b_{i+1}}}\frac{\partial \omega_{b_{i+1}}}{\partial l_{b_i}},\qquad i=1,\ldots,4.
\]

In Fig.~\ref{2-j kusok bipiramidy}
a fragment of bipyramid is presented that contains two identified
copies of each of edges $b_i$ and~$b_{i+1}$. One can see that
here a triple ``2nd case'' from
Section~\ref{sec obsch-proizv} takes place (which, of course, yields
again the 5th case).
Namely, the r\^ole of face $ABC$ from Fig.~\ref{fig peresek rebra}
is played in Fig.~\ref{2-j kusok bipiramidy} by each of faces
$C_0B_iB_{i+1}$, $B_{i+1}C_1C_0$ and~$C_1B_{i+1}B_{i+2}$. This gives
\begin{gather}
\left(F|_{\cal C}\right)_{i,i+1}=\frac{1}{6}\left(
\frac{V_{C_1B_iB_{i+1}B_{i+2}}}
{V_{C_0C_1B_{i+1}B_i} V_{C_0B_iB_{i+1}B_{i+2}}}
\right. \nonumber\\
\phantom{\left(F|_{\cal C}\right)_{i,i+1}=} \left.{}+ \frac{V_{C_0C_1B_{i+2}B_i}}
{V_{C_0C_1B_{i+1}B_i} V_{C_0C_1B_{i+2}B_{i+1}}} +
\frac{V_{C_0B_iB_{i+1}B_{i+2}}}{V_{C_0C_1B_{i+2}B_{i+1}}
V_{C_1B_iB_{i+1}B_{i+2}}} \right).
\label{eq pok2}
\end{gather}

Using the fact that points
$B_j$ and $C_j$ are located in two parallel planes,
we can again simplify the expression (\ref{eq pok2}) to
\begin{equation}
\left(F|_{\cal C}\right)_{i,i+1}=\frac{2}{6V_i} +\frac{2}{6V_{i+1}}.
\end{equation}

Perhaps these two types of matrix elements can be considered as
most complicated. Now we write out {\em all\/} matrix elements
of~$F|_{\cal C}$.

\subsubsection*{Submatrix {\mathversion{bold}$F|_{\cal C}$} for
{\mathversion{bold}$L(7,1)$} equals {\mathversion{bold}$1/6\, \times$}}
\begin{equation}
\left(
\begin{array}{@{}c@{\ }c@{\ }c@{\ }c@{\ }c@{}}
 {-} \frac{4}{V_{1}} {-} \frac{1}{V_{0}} {-} \frac{1}{V_{2}} &
 \frac{2}{V_{1}}  {+} \frac{2}{V_{2}} &
 {-} \frac{1}{V_{2}} & 0 & 0 \\[1ex]
 \frac{2}{V_{1}} {+} \frac{2}{V_{2}} &
 {-} \frac{4}{V_{2}} {-} \frac{1}{V_{1}} {-} \frac{1}{V_{3}} &
 \frac{2}{V_{2}} {+} \frac{2}{V_{3}} &
 {-} \frac{1}{V_{3}} & 0 \\[1ex]
 {-}\frac{ 1}{V_{2}} & \frac{2}{V_{2}} {+} \frac{2}{V_{3}} &
 {-} \frac{4}{V_{3}} {-} \frac{1}{V_{2}} {-} \frac{1}{V_{4}} &
 \frac{2}{V_{3}} {+} \frac{2}{V_{4}} &
 {-} \frac{1}{V_{4}}  \\[1ex]
0 & {-} \frac{1}{V_{3}} &
 \frac{2}{V_{3}} {+} \frac{2}{V_{4}} &
 {-} \frac{4}{V_{4}} {-} \frac{1}{V_{3}} {-} \frac{1}{V_{5}} &
 \frac{2}{V_{4}} {+} \frac{2}{V_{5}}  \\[1ex]
0 & 0 & {-} \frac{1}{V_{4}} &
 \frac{2}{V_{4}} {+} \frac{2}{V_{5}} &
 {-} \frac{4}{V_{5}} {-} \frac{1}{V_{4}} {-} \frac{1}{V_{6}}
\end{array}
\right).
\label{eq minor L(7,1)}
\end{equation}

\subsubsection*{Submatrix {\mathversion{bold}$F|_{\cal C}$} for
{\mathversion{bold}$L(7,2)$} equals {\mathversion{bold}$1/6\, \times$}}
\begin{equation}
\left(
\begin{array}{@{}c@{\ }c@{\ }c@{\ }c@{\ }c@{}}
 {-} \Sigma_{0123} &
 {-} \frac {1}{V_{2}} {+} \frac {1}{V_{1}} {+} \frac {1}{V_{3}} &
 \frac {1}{V_{2}} {+} \frac {1}{V_{3}} &
 {-} \frac {1}{V_{3}} &
 {-} \frac {1}{V_{0}} \\[1ex]
 {-} \frac {1}{V_{2}} {+} \frac {1}{V_{1}} {+} \frac {1}{V_{3}} &
 {-} \Sigma_{1234} &
 {-} \frac {1}{V_{3}} {+} \frac {1}{V_{2}} {+} \frac {1}{V_{4}} &
 \frac {1}{V_{3}} {+} \frac {1}{V_{4}} &
 {-} \frac {1}{V_{4}} \\[1ex]
 \frac {1}{V_{2}} {+} \frac {1}{V_{3}} &
 {-} \frac {1}{V_{3}} {+} \frac {1}{V_{2}} {+} \frac {1}{V_{4}} &
 {-} \Sigma_{2345} &
 {-} \frac {1}{V_{4}} {+} \frac {1}{V_{3}} {+} \frac {1}{V_{5}} &
 \frac {1}{V_{4}} {+} \frac {1}{V_{5}} \\[1ex]
 {-} \frac {1}{V_{3}} &
 \frac {1}{V_{3}} {+} \frac {1}{V_{4}} &
 {-} \frac {1}{V_{4}} {+} \frac {1}{V_{3}} {+} \frac {1}{V_{5}} &
 {-} \Sigma_{3456} &
 {-} \frac {1}{V_{5}} {+} \frac {1}{V_{4}} {+} \frac {1}{V_{6}} \\[1ex]
 {-} \frac {1}{V_{0}} &
 {-} \frac {1}{V_{4}} &
 \frac{1}{V_{4}} {+} \frac {1}{V_{5}} &
 {-} \frac {1}{V_{5}} {+} \frac {1}{V_{4}} {+} \frac {1}{V_{6}} &
 {-} \Sigma_{4560}
\end{array}
\right).
\label{eq minor L(7,2)}
\end{equation}
Here we had to introduce a brief notation
\[
\Sigma_{ijkl} \stackrel{\rm def}{=} \frac{1}{V_i} + \frac{1}{V_j} +
\frac{1}{V_k} + \frac{1}{V_l}
\]
in order that the matrix could be placed on the page.

Now we will be busy with the differential form (\ref{eq new30}).
{\em Before\/} calculating it we state the following

\medskip

\noindent
{\bf Proposition.} {\it The differential form (\ref{eq new30})
equals
\begin{equation}
|{\rm const}\cdot \rho\, d\rho \wedge \sigma \, d\sigma \wedge ds \wedge
d\alpha |
\label{eq *_20}
\end{equation}
(with notations taken from (\ref{eq *_10})).}

\medskip

\noindent
{\bf Proof.} First, we note that the rank of form (\ref{eq new30})
is equal to ${\mathop{\rm card}\nolimits}\, \overline{\cal C}$, that is $4$ in our case. Thus,
if expressed through $\rho$, $\sigma$, $s$ and $\alpha$, it
inevitably takes the form
\[
|f(\rho, \sigma, s, \alpha)\,
d\rho \wedge d\sigma \wedge ds \wedge d\alpha|.
\]
The form of function $f$ is easily obtained from reasoning of the same type
as in Section~5 of paper~\cite{1}: differential form~(\ref{eq new30}) must be
represented as the exterior product of the volume form for point~$C$, i.e.\
$dx_C \wedge dy_C \wedge dz_C = \sigma\, d\sigma \wedge ds \wedge d\alpha $
by a factor not depending on the location of point~$C$ with respect to
point~$B$ and symmetry axis~$z$. This must also hold for
point~$B$. Thus, we are led unambiguously to formula~(\ref{eq *_20})
for the form~(\ref{eq new30}). The proposition is proved.
\hfill \qed

\medskip

In the case we are considering,
it is natural to take the constant entering in~(\ref{eq *_20}) as
our invariant (which, we recall, we put in correspondence to the {\em pair\/}
$(M,\varphi)$, where $\varphi$ is a homomorphism from
$\pi_1(M)$ to~$E_3$).

In the numerator of formula~(\ref{eq new30}) we see the exterior product
$\displaystyle\bigwedge_{\overline{\cal C}} l\, dl$ which in our case,
as one can check, has the form
(we write it out for arbitrary $p$ and~$q$)
\begin{gather*}
 l_a\, dl_a \wedge l_{b_0}\, dl_{b_0} \wedge l_{b_{p-1}}\, dl_{b_{p-1}}
\wedge l_c\, dl_c \\
\qquad {} = 8 R \sin^2 \frac{\pi k}{p}
\sin \frac{\pi qk}{p}
\cos \left(\alpha + \frac{\pi k}{p}\right)
\rho\, d\rho \wedge \sigma \, d\sigma \wedge ds \wedge d\alpha
\end{gather*}
($R$ was introduced in formula (\ref{eq *_R})). This shows that the
{\bf dependence of the denominator of formula~(\ref{eq new30}) on~$\alpha$
must, too, reduce to a factor~$\cos (\alpha + \pi k/p)$}, in order that
this factor could cancel out and we get a formula of type~(\ref{eq *_20}).
We hope that, with this remark, our calculations of the
expression~(\ref{eq new30}), for which we actively used the
Maple program (the most difficult step was the calculation of determinants
of matrices (\ref{eq minor L(7,1)}) and~(\ref{eq minor L(7,2)})),
will be repeatable.

Here are the results for our invariant, i.e.\ value ``$\rm const$'' from
formula~(\ref{eq *_20}), which we will now call $I_k$ because it
depends on the integer $k$ from formulas~(\ref{eq *_10}) that determines
the homomorphism~$\varphi$. We recall that it is enough to consider
$k=1,2$ or~$3$ for $L(7,q)$. So, with the precision of ten
significant digits,
\begin{gather}
I_1 \bigl( L(7,1) \bigr)  =  0.08100567416,
\label{eq I_1 L(7,1)} \\
I_2 \bigl( L(7,1) \bigr)  =  0.8540328192, \\
I_3 \bigl( L(7,1) \bigr)  =  2.064961508, \\
I_1 \bigl( L(7,2) \bigr)  =  0.2630237713, \\
I_2 \bigl( L(7,2) \bigr)  =  1.327985278, \\
I_3 \bigl( L(7,2) \bigr)  =  0.4089909518.
\label{eq I_3 L(7,2)}
\end{gather}

It is obvious that $L(7,1)$ cannot be homeomorphic to~$L(7,2)$.
Exact expressions for values of $I_k$ in terms of trigonometric functions
are presented in Discussion, formula~(\ref{eq I_k obsch}),
together with a conjecture about their form for an arbitrary lens space.

\section{Discussion}
\label{sec discussion}

Invariants of the stated type can be calculated for any 3-dimensional
closed orientable manifolds. It looks very plausible that they can be
also generalized to manifolds with boundary and not necessarily
orientable. It looks clear that the differential form~(\ref{eq *_20})
whose very appearance is connected with the presence
of one symmetry axis in the image in~$\mathbb R^3$ of the universal covering
of manifold pre-triangulation will have to be replaced with
other forms, as soon as the fundamental group of the
manifold is different from~$\mathbb Z_p$ (recall also
the differential form \cite[(33)]{1} that corresponded to the
trivial homomorphism~$\varphi$).

As for the lens spaces with arbitrary $p$ and $q$, we have the
following conjecture about the values of invariant $I_k$ for them:
\begin{equation}
I_k\bigl( L(p,q) \bigr) = \frac{16}{p} \sin^2 \frac{\pi k}{p}\,
\sin^2 \frac{\pi qk}{p},\qquad k=1, \ldots, \left[ \frac{p}{2} \right].
\label{eq I_k obsch}
\end{equation}
Approximate values of exactly these numbers were presented in formulas
(\ref{eq I_1 L(7,1)})--(\ref{eq I_3 L(7,2)}).
Note that formula~(\ref{eq I_k obsch}) agrees with the known
fact~\cite{HW,milnor} that lens spaces $L(p,q_1)$ and $L(p,q_2)$
are homeomorphic if (and only if)
\[
q_2 = \pm q_1^{\pm 1} \pmod{p}.
\]
Indeed, if $q_2 = \pm q_1\pmod{p}$, then, obviously,
$I_k\bigl( L(p,q_1) \bigr) = I_k\bigl( L(p,q_2) \bigr)$.
If, on the other hand,
$q_2 = \pm q_1^{-1} \pmod{p}$, then
$I_k\bigl( L(p,q_1) \bigr) = I_{kq_1} \bigl( L(p,q_2) \bigr)$,
i.e.\ the set of numbers corresponding to the space is the same.

Our invariants may be interesting also because they depend not only
on a manifold~$M$ but also on a homomorphism
$\varphi\colon\; \pi_1(M)\to E_3$.
In particular, they make it possible to distinguish the generators
of fundamental groups for lens spaces.

Now we mention the following unsolved and intriguing problems.
\vspace{-2mm}
\begin{itemize}
\itemsep0mm
\item
In preprint~\cite{korepanov 4s} a ``higher analogue of pentagon
equation'' is presented~--- an algebraic relation corresponding
to a $3\to 3$ re-building of a cluster of 4-simplices (expressed, again,
in terms of Euclidean metric values assigned to those simplices).
It would be very interesting to construct a 4-manifold invariant on this
basis. If our classical invariants work for higher-dimensional manifolds,
which looks plausible, this would be an advantage with respect to the known
quantum invariants because, for these latter,
it is quite unclear how to
generalize them nontrivially to higher dimensions.
Maybe, equations of pentagon type that arise in connection with
multidimensional topology will give an idea of how to construct
integrable models in multidimensional {\em mathematical physics\/}
as well.

\item
In paper~\cite{SL_2}, we have constructed one more classical solution
of pentagon equation which we called $SL(2)$-solution. We must recognize
that we were not yet able to ``globalize'' it, that is build on its
basis invariants of 3-manifolds, not to say of passing to higher dimensions.

\item
The third intriguing question is whether we can invent some quantum objects
of which our classical solutions (presented in \cite{1,korepanov 4s}
and~\cite{SL_2}) are quasiclassical limits.
\vspace{-2mm}
\end{itemize}

\subsection*{Acknowledgements}
One of the authors (I~Korepanov) acknowledges
the partial financial support from Russian Foundation for
Basic Research under Grant no.~01-01-00059.

\subsection*{Note added in proof}

For a given $p$, formula~(\ref{eq I_k obsch}) gives, essentially,
the squared {\em Reidemeister torsion\/} of $L(p,q)$ multiplied
by some constant. We thank Yuka Taylor for this intriguing
observation.

\label{korepanov-lastpage}

\end{document}